\pdfoutput=1
\RequirePackage{ifpdf}
\ifpdf 
\documentclass[pdftex]{sigma}
\else
\documentclass{sigma}
\fi

\numberwithin{equation}{section}

\newtheorem{Theorem}{Theorem}[section]
\newtheorem*{Theorem*}{Theorem}

\newtheorem{Lemma}[Theorem]{Lemma}
\newtheorem{Conjecture}[Theorem]{Conjecture}
\newtheorem{Proposition}[Theorem]{Proposition}
 { \theoremstyle{definition}

 }
\usepackage{tikz}
\usetikzlibrary{arrows,patterns,decorations.pathmorphing,backgrounds,positioning,fit,petri,calligraphy}
\usetikzlibrary {arrows.meta}

\begin{document}
\allowdisplaybreaks

\newcommand{\arXivNumber}{2101.11708}

\renewcommand{\PaperNumber}{002}

\FirstPageHeading

\ShortArticleName{A Cable Knot and BPS-Series}

\ArticleName{A Cable Knot and BPS-Series}

\Author{John CHAE}

\AuthorNameForHeading{J.~Chae}

\Address{Department of Mathematics, Univeristy of California Davis, Davis, USA}
\Email{\href{mailto:yjchae@ucdavis.edu}{yjchae@ucdavis.edu}}

\ArticleDates{Received August 03, 2022, in final form January 05, 2023; Published online January 13, 2023}

\Abstract{A series invariant of a complement of a knot was introduced recently. The invariant for several prime knots up to ten crossings have been explicitly computed. We present the first example of a satellite knot, namely, a cable of the figure eight knot, which has more than ten crossings. This cable knot result provides nontrivial evidence for the conjectures for the series invariant and demonstrates the robustness of integrality of the quantum invariant under the cabling operation. Furthermore, we observe a relation between the series invariant of the cable knot and the series invariant of the figure eight knot. This relation provides an alternative simple method of finding the former series invariant.}

\Keywords{knot complement; quantum invariant; $q$-series; Chern--Simons theory; categorification}

\Classification{57K10; 57K16; 57K31; 81R50}

\section{Introduction}

Inspired by a categorification of the Witten--Reshitikhin--Turaev invariant of a closed oriented 3-manifold~\cite{RT2,RT1, W} in \cite{GPPV, GPV}, a two variable series invariant $F_K(x,q)$ for a complement of a~knot $M^{3}_K$ was introduced in \cite{GM}. Although its rigorous definition is yet to be found, it possesses various properties such as the Dehn surgery formula and the gluing formula. This knot invariant $F_K$ takes the form\footnote{Implicitly, there is a choice of group; originally, the group used is ${\rm SU}(2)$.}
\begin{gather}\label{eq1.1.1}
F_K(x,q)= \frac{1}{2} \sum_{\substack{m \geq 1 \\ m \ \text{odd}}}^{\infty} \big(x^{m/2}-x^{-m/2}\big)f_{m}(q) \in \frac{1}{2^{c}} q^{\Delta} \mathbb{Z}\big[x^{\pm 1/2}\big]\big[\big[q^{\pm 1}\big]\big],
\end{gather}
where $f_{m}(q)$ are Laurent series with integer coefficients,\footnote{They can be polynomials for monic Alexander polynomial of $K$, see Section~\ref{section3.2}.} $c \in \mathbb{Z}_{+}$ and $\Delta \in \mathbb{Q}$. Moreover, $x$-variable is associated to the relative $\operatorname{Spin}^c \big(M^{3}_K, T^2\big)$-structures, which is affinely isomorphic to $H^2\big(M^{3}_K, T^2 ; \mathbb{Z}\big) \cong H_1\big(M^{3}_K;\mathbb{Z}\big) \cong \mathbb{Z}$. It is infinite cyclic, which is reflected as a series in $F_K$. The rational constant $\Delta$ was investigated in \cite{GPP}, which elucidated its intimate connection to the d-invariant (or the correction term) in certain versions of the Heegaard Floer homology $\big(HF^{\pm}\big)$ for rational homology spheres. The physical interpretation of the integer coefficients in $f_{m}(q)$ are number of BPS states of 3d $\mathcal{N}=2$ supersymmetric quantum field theory on $M^{3}_{K}$ together with boundary conditions on $\partial M^{3}_{K}$. Furthermore, it was conjectured that $F_K$ also satisfies the Melvin--Morton--Rozansky conjecture~\cite{MM,R,R2} (proven in~\cite{BG}):
\begin{Conjecture}[{\cite[Conjecture 1.5]{GM}}]\label{conjecture3} For a knot $K \subset$ $S^3$, the asymptotic expansion of the knot invariant $F_{K}\big(x,q={\rm e}^{\hbar}\big)$ about $\hbar =0$ coincides with the Melvin--Morton--Rozansky $($MMR$)$ expansion of the colored Jones polynomial in the large color limit:
\begin{gather}\label{eq1.1.2}
\frac{F_{K}\big(x,q={\rm e}^{\hbar}\big)}{x^{1/2}-x^{-1/2}} = \sum_{r=0}^{\infty} \frac{P_{r}(x)}{\Delta_K(x)^{2r+1}}\hbar^r,
\end{gather}
where $x={\rm e}^{n\hbar}$ is fixed, n is the color of $K$, $P_{r}(x) \in \mathbb{Q} \big[x^{\pm 1}\big]$, $P_{0}(x)=1$ and $\Delta_K (x)$ is the $($symmetrized$)$ Alexander polynomial of~$K$.
\end{Conjecture}

Additionally, motivated by the quantum volume conjecture/AJ-conjecture~\cite{G,Gukov2} (explained in Section~\ref{sec22}), it was conjectured that $F_K$-series is $q$-holonomic:
\begin{Conjecture}[{\cite[Conjecture 1.6]{GM}}] For any knot $K \subset$ $S^3$, the normalized series $f_{K}(x,q)$ satisfies a linear recursion relation generated by the quantum A-polynomial of $K$ $\hat{A}_K(q,\hat{x},\hat{y})$:
\begin{gather}\label{eq1.1.3}
\hat{A}_{K}(q, \hat{x},\hat{y}) f_{K}(x,q) = 0,
\end{gather}
where $f_{K}:=F_{K}(x,q)/\big(x^{1/2}-x^{-1/2}\big)$.
\end{Conjecture}

The actions of $\hat{x}$ and $\hat{y}$ are
\begin{equation*}
\hat{x} f_{K}(x,q)= x f_{K}(x,q), \qquad \hat{y}f_{K}(x,q)= f_{K}(xq,q).
\end{equation*}
$F_K$-series has been computed for several prime knots up to ten crossings in \cite{GHNPPS, GM,K,P2}. They include the torus knots, the figure eight knot in \cite{GM,P,P2}. Positive braid knots ($10_{139}$, $10_{152}$), strongly quasipositive braids knots ($m(10_{145})$, $10_{154}$, $10_{161}$), double twist knots ($m(5_{2})$, $m(7_{3})$, $m(7_{4})$), and a few more prime knots ($m(7_{5})$, $m(8_{15})$) were examined in~\cite{P2}. Furthermore, the series for $5_{2}$ and $6_{2}$ were calculated in~\cite{EGGKPSS}.

In this paper, we verify the above conjectures by computing the $F_K$-series for $(9,2)$-cabling of the figure eight knot and we compare our result to that of the figure eight knot. Furthermore, we conjecture the form of the $F_K$-series for a family of a cable knot of the figure eight.

The rest of the paper is organized as follows. In Section~\ref{section2}, we review the satellite operation on a knot and the recursion ideal of the quantum torus. In Section~\ref{sec3}, we analyze knot polynomials of the cable knot of the figure eight. In Section~\ref{sec4}, we derive the recursion relation for the cable knot. Then we deduce $\hbar$ expansion from the recursion in Section~\ref{sec5}. In Section~\ref{sec6}, consequences of the cabling operation are discussed and we propose a conjecture about a family of a cable knot. Finally, in Section~\ref{sec7}, we state a relation between the series invariant of the cable knot and the series invariant of the figure eight knot and conjecture about other cabling of the figure eight knot.

\section{Background}\label{section2}

\subsection{Satellites}

The satellite operation consists of a pattern knot $P$ in the interior of the solid torus $S^1 \times D^2$, a~companion knot $K^{\prime}$ in the $S^3$ and an canonical identification $h_{K^{\prime}}$
\begin{equation}\label{eq2.2.1}
h_{K^{\prime}} \colon \ S^1 \times D^2 \longrightarrow \nu(K^{\prime}) \subset S^3,
\end{equation}
where $\nu(K^{\prime})$ is the tubular neighborhood of $K^{\prime}$.
A well-known example of satellite knots is a~cable knot $h_{K^{\prime}}(P)=C_{(r,s)}(K^{\prime})$ that is obtained by choosing P to be the $(r,s)$-torus knot pushed into the interior of the $S^1 \times D^2$. This map $h_{K^{\prime}}$ has been investigated in \cite{Levine, Miller, MP}.
\begin{figure}[t]
 \centering
		\includegraphics[scale=0.4]{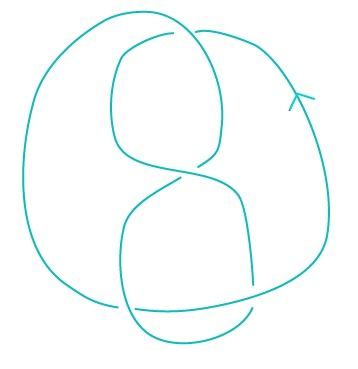}
		\qquad
 \includegraphics[scale=0.3]{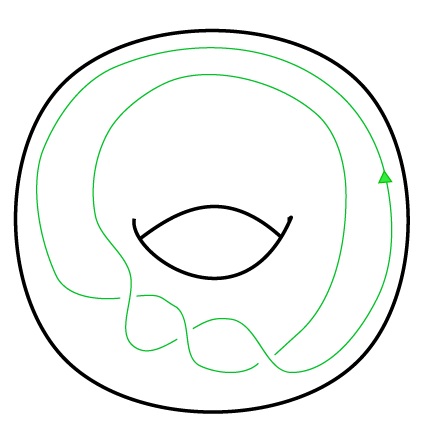}
 \qquad
		\includegraphics[scale=0.2]{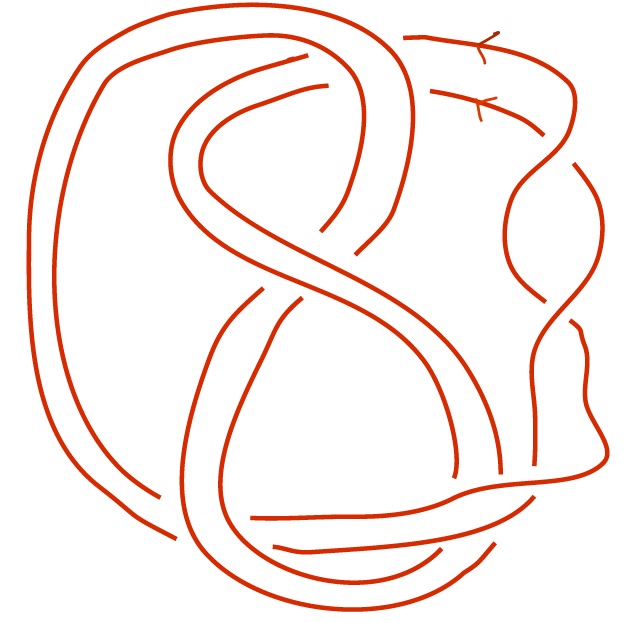}
 \caption{A companion $K^{\prime}$ (left), pattern knot $P$ (center) and satellite knot $P(K^{\prime})$ (right).}
\end{figure}

\subsection{Quantum torus and recursion ideal}\label{sec22}

Let $\mathcal{T}$ be a quantum torus
\begin{equation*}
\mathcal{T} : = \mathbb{C} \big[ t^{\pm 1} \big] \big\langle M^{\pm 1}, L^{\pm 1}\big\rangle / \big(LM- t^2 ML\big).
\end{equation*}
The generators of the noncommutative ring $\mathcal{T}$ acts on a set of discrete functions, which are colored Jones polynomials $J_{K,n} \in \mathbb{Z}\big[t^{\pm 1}\big]$ in our context, as
\begin{equation*}
M J_{K,n}= t^{2n} J_{K,n}, \qquad L J_{K,n}= J_{K,n+1}.
\end{equation*}
The recursion (annihilator) ideal $\mathcal{A}_{K}$ of $J_{K,n}$ is the left ideal $\mathcal{A}_{K}$ in $\mathcal{T}$ consisting of operators that annihilates $J_{K,n}$:
\begin{equation*}
\mathcal{A}_{J_{K,n}} : = \big\{ \alpha_{K} \in \mathcal{T} \, |\, \alpha_{K} J_{K,n} = 0 \big\}.
\end{equation*}
It turns out that $\mathcal{A}_{K}$ is not a principal ideal in general. However, by adding inverse of polynomials of $t$ and $M$ to $\mathcal{T}$~\cite{G},
we obtain a principal ideal domain $\tilde{\mathcal{T}}$
\begin{equation*}
\tilde{\mathcal{T}} : = \bigg\{ \sum_{j \in \mathbb{Z}} a_{j}(M) L^{j} \,\big|\, a_{j}(M) \in \mathbb{C}\big[t^{\pm 1}\big](M),\, a_{j}= \text{almost always} \ 0 \bigg\}
\end{equation*}
Using $\tilde{\mathcal{T}}$ we get a principal ideal $\tilde{\mathcal{A}_{K}}:= \tilde{\mathcal{T}}\mathcal{A}_{K}$ generated by a single polynomial $\hat{A}_{K}$
\begin{equation*}
\hat{A}_{K}(t,M,L)= \sum_{j=0}^{d} a_{j}(t,M)L^{j}.
\end{equation*}
This $\hat{A}_{K}$ polynomial is a noncommutative deformation of a classical A-polynomial of a knot~\cite{CCGLS} (see also \cite{CL}). Alternative approaches to obtain $\hat{A}_K(t,M,L)$ are by quantizing the classical A-polynomial curve using a twisted Alexander polynomial or applying the topological recursion~\cite{GS}. A conjecture called AJ conjecture/quantum volume conjecture was proposed in \cite{G,Gukov2} via different approaches:
\begin{Conjecture}
For any knot $K \subset$ $S^3$, $\hat{A}_{K}$ $(t=-1, L, M)$ reduces to the $($classical$)$ A-poly\-no\-mi\-al curve $A_{K}(L,M)$ up to a solely $M$-dependent overall factor.
\end{Conjecture}

In other words, $J_{K,n}(t)$ satisfies a linear recursion relation generated by $\hat{A}_{K}(t,M,L)$. This property of $J_{K,n}$ is often called $q$-holonomic~\cite{GL}. The conjecture was confirmed for a variety of knots~\cite{DGLZ,G,GK,GS2,H,LT,RZ,Tran2}.

\section{Knot polynomials}\label{sec3}

In this section we will analyze the colored Jones polynomial and the Alexander polynomial of a cable knot to show that the former satisfies the MMR expansion and the latter is monic. Furthermore, the MMR expansion enables us to read off the initial condition that is needed in Section~\ref{sec5}.

\subsection{The colored Jones polynomial}

For $(r,2)$-cabling of the figure eight knot $4_1$, we set $P=T(r,2)$ and $K^{\prime}=4_1$ in \eqref{eq2.2.1}. The cabling formula for an unnormalized $\mathfrak{sl}_2(\mathbb{C})$ colored Jones polynomial of a $(r,2)$-cabling of a~knot~$K^{\prime}$ in~$S^3$ is~\cite{Tran}
\begin{equation*}
\tilde{J}_{C_{(r,2)}(K^{\prime}),n}(q)= q^{-\frac{r}{2} (n^2-1 )} \sum_{w=1}^{n} (-1)^{r (n-w)} q^{\frac{r}{2} w(w-1)} \tilde{J}_{K^{\prime}, (2w-1)}(q),\qquad |r| > 8\ \text{and odd}.
\end{equation*}

\begin{figure}[h!]\vspace{-5mm}
\centering
 \begin{tikzpicture}[scale=1,smooth]
 \node at (0,0) {\includegraphics[scale=0.4]{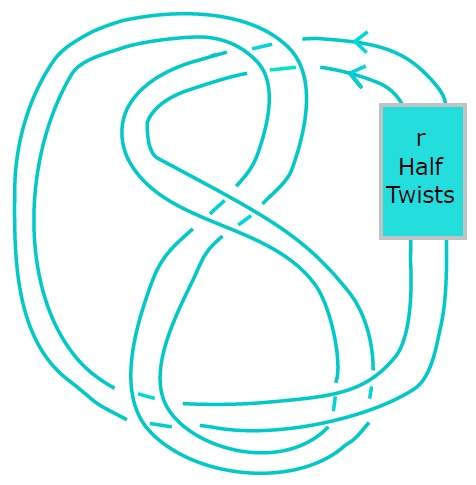}};
\draw [cyan!60,fill=cyan!60] (1.45,0.0) rectangle (2.45,1.5);
\node at (1.95,0.75) {\small $\begin{array}{c} r \\ \text{half} \\ \text{twists}\end{array}$};

\end{tikzpicture}
\caption{$(r,2)$-cable of the figure eight knot. }
\end{figure}

Its application to $K=C_{(9,2)}(4_1)$,\footnote{This cabling parameters correspond to $9_1$ for the pattern knot. We assume 0-framing for $4_1$.} whose diagram has 25 crossings, is
\begin{gather*}
		 \tilde{J}_{K,n}(q) = q^{-\frac{9}{2} (n^2-1)} \sum_{w=1}^{n} \left[\vphantom{\sum_{r=0}^{2w-2}} (-1)^{(n-w)} q^{\frac{9}{2} w(w-1)} [2w-1] \right.\\
\left.\hphantom{\tilde{J}_{K,n}(q) = q^{-\frac{9}{2} (n^2-1)}}{}
\times \sum_{r=0}^{2w-2} \prod_{k=1}^{r} \bigl( -q^{-k}-q^k+q^{1-2 w}+q^{2 w-1} \bigr) \right].
\end{gather*}
Using the (0-framed) unknot $U$ value
\begin{equation*}
\tilde{J}_{U,n}(t) = \frac{t^{2n} - t^{-2n}}{t^2 - t^{-2}},
\end{equation*}
together with $q=t^4$, the first few normalized polynomials $J_{K,n}(q)$ can be written as
\begin{gather*}	
		 J_{K,1}(q)=1,\\
		 J_{K,2}(q)=q^2-q+\frac{1}{q^4}+\frac{1}{q^6}-\frac{1}{q^7}+\frac{1}{q^8}-\frac{1}{q^9}+\frac{1}{q^{12}}-\frac{1}{q^{13}},\\
J_{K,3}(q) = q^{12}-q^{11}-q^{10}+q^9 -q^8 +q^7+q^6-q^5+q^2-1+\frac{1}{q^8}+\frac{1}{q^{11}}-\frac{1}{q^{13}}+\frac{1}{q^{14}}-\frac{1}{q^{16}}\\
\hphantom{J_{K,3}(q) =}{}
+\frac{1}{q^{17}} -\frac{1}{q^{18}}-\frac{1}{q^{19}} +\frac{2}{q^{20}} -\frac{1}{q^{21}}+\frac{1}{q^{23}}-\frac{1}{q^{24}}+\frac{1}{q^{25}}+\frac{1}{q^{26}}-\frac{2}{q^{27}}-\frac{1}{q^{28}}+\frac{1}{q^{29}}\\
\hphantom{J_{K,3}(q) =}{}-\frac{1}{q^{30}} +\frac{2}{q^{32}}-\frac{1}{q^{33}} -\frac{1}{q^{34}}+\frac{1}{q^{35}}.
\end{gather*}
\begin{Proposition}
The $\hbar$ expansion of the above $J_{K,n}(q)$ is given by
\begin{gather}		
	J_{K,n}\big({\rm e}^{\hbar}\big) = 1 + \big( 6 - 6 n^2 \big) \hbar^2 + \big( {-}42 + 42 n^2 \big) \hbar^3 + \left( \frac{801}{2} - 462 n^2 + \frac{123}{2} n^4 \right) \hbar^4\nonumber\\
\hphantom{J_{K,n}\big({\rm e}^{\hbar}\big) =}{}
	 + \left( - \frac{8451}{2} + 5173 n^2 - \frac{1895}{2} n^4 \right) \hbar^5 + \left( \frac{3111491}{60} - \frac{132779}{2} n^2 + 14986 n^4 \right.\nonumber\\
\left. \hphantom{J_{K,n}\big({\rm e}^{\hbar}\big) =}{} - \frac{27281}{60} n^6 \right) \hbar^6
	 + \left( -\frac{14631401}{20} + \frac{19399417}{20} n^2 - \frac{3028829}{12} n^4 + \frac{840097}{60} n^6 \right) \hbar^7\nonumber\\
\hphantom{J_{K,n}\big({\rm e}^{\hbar}\big) =}{}
	 + \left( \frac{39069313501}{3360} - \frac{950122877}{60} n^2 + \frac{54585517}{12} n^4 - \frac{1725671}{5} n^6\right.\nonumber\\
\left. \hphantom{J_{K,n}\big({\rm e}^{\hbar}\big) =}{}
 + \frac{13273763}{3360} n^8 \right) \hbar^8 + \cdots.\label{eq3.3.1}
\end{gather}
\end{Proposition}

We see that, at each $\hbar$ order, the degree of the polynomial in $n$ is at most the order of $\hbar$, which is an equivalent characterization of the MMR expansion of the colored Jones polynomial of a~knot. Secondly, as a consequence of the cabling, odd powers of $\hbar$ appear in the expansion, even though they are absent in the case of the figure eight knot~\cite{GM}.

\subsection{The Alexander polynomial}\label{section3.2}

The cabling formula for the Alexander polynomial of a knot $K$ is~\cite{Hedden}
\begin{equation*}
\Delta_{C_{(p,q)}(K)}(t)= \Delta_K (t^p) \Delta_{T_{(p,q)}}(t), \qquad 2 \leq p < |q|, \quad \gcd (p,q)=1,
\end{equation*}
where $\Delta(t)$ is the symmetrized Alexander polynomial and $T_{(p,q)}$ is the $(p,q)$ torus knot. Note that our convention for the parameters of the torus knot are switched (i.e., $p\equiv 2$, $q\equiv r$).
\begin{Lemma}
The symmetrized Alexander polynomial of $C_{(9,2)}(4_1)$ is as follows:
\begin{align*}		
\Delta_{C_{(9,2)}(4_1)}(x) & = \Delta_{4_1} \big(x^2\big) \Delta_{T_{(2,9)}}(x)\\
& = -x^6-\frac{1}{x^6}+x^5+\frac{1}{x^5}+2 x^4+\frac{2}{x^4}-2 x^3-\frac{2}{x^3}+x^2+\frac{1}{x^2}-x-\frac{1}{x} +1.
\end{align*}
\end{Lemma}

From this Alexander polynomial its symmetric expansion about $x=0$ (in $x$) and $x=\infty$ (in~$1/x$) in the limit of $\hbar \rightarrow 0$ can be computed:
\begin{gather}		
\lim_{q \rightarrow 1} 2F_{K}(x,q) = 2\operatorname{s.e.} \left( \frac{x^{1/2}-x^{-1/2}}{\Delta_{K}(x)} \right)\nonumber\\
 \hphantom{\lim_{q \rightarrow 1} 2F_{K}(x,q)}{}
 = x^{11/2} -\frac{1}{x^{11/2}} +2 x^{15/2} - \frac{2}{x^{15/2}} + 5 x^{19/2} -\frac{5}{x^{19/2}} +13 x^{23/2} - \frac{13}{x^{23/2}}\nonumber\\
\hphantom{\lim_{q \rightarrow 1} 2F_{K}(x,q)=}{} + 34 x^{27/2} - \frac{34}{x^{27/2}} -x^{29/2} + \frac{1}{x^{29/2}} +89 x^{31/2} - \frac{89}{x^{31/2}} -2 x^{33/2} + \frac{2}{x^{33/2}}\nonumber\\
\hphantom{\lim_{q \rightarrow 1} 2F_{K}(x,q)=}{} + 233 x^{35/2} - \frac{233}{x^{35/2}} -5 x^{37/2} + \frac{5}{x^{37/2}} +610 x^{39/2}- \frac{610}{x^{39/2}} + \cdots\nonumber\\
\hphantom{\lim_{q \rightarrow 1} 2F_{K}(x,q)=}{}
 \in \mathbb{Z} \big[ \big[ x^{\pm 1/2} \big] \big].\label{eq3.3.2}
\end{gather}
The coefficients in the expansions are integers and hence the Alexander polynomial is monic, which is a necessary condition for $f_m(q)$'s in \eqref{eq1.1.1} to be polynomials.

\section{The recursion relation}\label{sec4}
The quantum (or noncommutative) A-polynomial of a class of cable knot $C_{(r,2)}(4_1)$ in $S^3$ having minimal $L$-degree is given by~\cite{Ruppe}
\begin{equation}
\hat{A}_{K} (t,M,L)=(L-1)B(t,M)^{-1}Q(t,M,L) \big( M^r L + t^{-2r}M^{-r} \big) \in \tilde{\mathcal{A}}_{K},\label{eq4.4.1}
\end{equation}
where
\begin{gather*}
Q(t,M,L)= Q_2(t,M)L^2 + Q_1(t,M)L + Q_0(t,M),\qquad B(t,M):= \sum_{j=0}^{2}c_{j}b\big(t,t^{2j+2}M^2\big),
\\
b(t,M)=\frac{M\big(1+t^4 M^2\big)\big({-}1+t^4 M^4\big)\big({-}t^2 + t^{14} M^4\big)}{t^2 - t^{-2}},
\\
c_{0}=\hat{P}_0 \big(t, t^4 M^2\big)\hat{P}_1 \big(t,t^6 M^2\big),\qquad c_{1}= - \hat{P}_1\big(t,t^2 M^2\big)\hat{P}_{1}\big(t,t^6 M^2\big),\\
 c_{2}=\hat{P}_1\big(t,t^2 M^2\big)\hat{P}_2\big(t,t^4 M^2\big).
\end{gather*}
The definitions of the operators $\hat{P}_i$ are written in Appendix~\ref{appendixA}. For $K=C_{(9,2)}(4_1)$, applying~\eqref{eq4.4.1} to $f_K(x,q)$ together with $x=q^n$ yields via \eqref{eq1.1.3}
\begin{gather}
\alpha(x,q)F_K(x,q) + \beta(x,q)F_K(xq,q)+ \gamma(x,q) F_K\big(xq^2,q\big) \nonumber\\
\qquad{}+ \delta(x,q)F_K\big(xq^3,q\big) + F_K\big(xq^4,q\big)=0,\label{eq4.4.2}
\end{gather}
where $\alpha$, $\beta$, $\gamma$, $\delta$ functions and their $\hbar$ series are documented in \cite{C}. From \eqref{eq4.4.2} we find the recursion relation for $f_{m}$.
\begin{Theorem}
The recursion relation for $f_{m}(q) \in \mathbb{Z}[q^{\pm 1}]$ of the above $F_K (x,q)$ is given by
\begin{gather}
f_{m+98} (q) = \frac{-1}{q^{\frac{109+m}{2}} \big(1-q^{\frac{87+m}{2}}\big)} \big[ t_{2} f_{m+94}+ t_{4} f_{m+90} + t_{6} f_{m+86} + t_{8} f_{m+82} + t_{9} f_{m+80}\nonumber\\
 \hphantom{f_{m+98} (q) =}{}
+ t_{10} f_{m+78} + t_{11} f_{m+76} + t_{12} f_{m+74} + t_{13} f_{m+72}+ t_{14} f_{m+70} + t_{15} f_{m+68}\nonumber\\
 \hphantom{f_{m+98} (q) =}{} + t_{16} f_{m+66} + t_{17} f_{m+64} + t_{18} f_{m+62} + t_{19} f_{m+60} + t_{20} f_{m+58} + t_{21} f_{m+56}\nonumber\\
 \hphantom{f_{m+98} (q) =}{} + t_{22} f_{m+54}+ t_{23} f_{m+52} + t_{24} f_{m+50} + t_{25} f_{m+48} + t_{26} f_{m+46} + t_{27} f_{m+44}\nonumber\\
 \hphantom{f_{m+98} (q) =}{}
 + t_{28} f_{m+42} + t_{29} f_{m+40} + t_{30} f_{m+38} + t_{31} f_{m+36} + t_{32} f_{m+34}+ t_{33} f_{m+32}\nonumber\\
 \hphantom{f_{m+98} (q) =}{}
+ t_{34} f_{m+30} + t_{35} f_{m+28} + t_{36} f_{m+26} + t_{37} f_{m+24}+ t_{38} f_{m+22} + t_{39} f_{m+20}\label{eq4.4.3}\\
\hphantom{f_{m+98} (q) =}{}+ t_{40} f_{m+18} + t_{41} f_{m+16}+ t_{43} f_{m+12} + t_{45} f_{m+8} + t_{47} f_{m+4} + t_{49} f_{m} \big] \in \mathbb{Z}\big[q^{\pm 1} \big], \nonumber
\end{gather}
where $t_{v}=t_{v}(q,q^m)$'s are listed in~{\rm \cite{C}}. The initial data for \eqref{eq4.4.3} were found using the $\hbar$-ex\-pan\-sions of \eqref{eq4.4.2}. An example of the expansion is written in Section~{\rm \ref{sec5}} for~$F_K (x,q)$. Using the recursion relation \eqref{eq4.4.3} and the initial data documented in~{\rm \cite{C}}, $F_{K}(x,q)$ can be obtained to any desired order in~$x$.
\end{Theorem}

\section{An expansion of a knot complement}\label{sec5}
We next compute a series expansion of the $F_K$ of complement of the cable knot $K$. Specifically, a~straightforward computation from \eqref{eq4.4.2} yields an ordinary differential equation (ODE) for~$P_m(x)$ at each $\hbar$ order. Using the initial conditions for the ODEs obtained from \eqref{eq3.3.1}
\begin{gather*}
P_{1}(1)=0,\qquad\! P_{2}(1)=6,\qquad\! P_{3}(1)= -42,\qquad\! P_{4}(1)= \frac{801}{2},\qquad\! P_{5}(1)= -\frac{8451}{2},\qquad\! \dots,
\end{gather*}
we find that
\begin{gather*}
P_{1}(x) = 5 x^{12}+\frac{5}{x^{12}}-10 x^{11}-\frac{10}{x^{11}}-13 x^{10}-\frac{13}{x^{10}}+36 x^9+\frac{36}{x^9}-10 x^8-\frac{10}{x^8}-16 x^7\\
\hphantom{P_{1}(x) =}{} -\frac{16}{x^7}+15x^6 +\frac{15}{x^6} -14 x^5-\frac{14}{x^5}+16 x^4+\frac{16}{x^4}-18 x^3-\frac{18}{x^3}+19 x^2+\frac{19}{x^2}\\
\hphantom{P_{1}(x) =}{}
-20 x-\frac{20}{x} +20,\\
P_{2}(x) = \frac{25 x^{24}}{2}+\frac{25}{2 x^{24}}-50 x^{23}-\frac{50}{x^{23}}-14 x^{22}-\frac{14}{x^{22}}+306 x^{21}+\frac{306}{x^{21}}-\frac{641x^{20}}{2}-\frac{641}{2 x^{20}}\\
\hphantom{P_{2}(x) =}{}
-448 x^{19} -\frac{448}{x^{19}} +\frac{2011 x^{18}}{2}+\frac{2011}{2 x^{18}}-358 x^{17}-\frac{358}{x^{17}}-522 x^{16} -\frac{522}{x^{16}}+612 x^{15}\\
\hphantom{P_{2}(x) =}{}
+\frac{612}{x^{15}}-\frac{589 x^{14}}{2} -\frac{589}{2 x^{14}} +508 x^{13} +\frac{508}{x^{13}} -\frac{3325 x^{12}}{2}-\frac{3325}{2 x^{12}}+1648 x^{11} +\frac{1648}{x^{11}}\\
\hphantom{P_{2}(x) =}{}
+1538 x^{10}+\frac{1538}{x^{10}}-3932 x^9
 -\frac{3932}{x^9} +1574 x^8 +\frac{1574}{x^8} +1670 x^7+\frac{1670}{x^7}-1798 x^6 \\
\hphantom{P_{2}(x) =}{}
 -\frac{1798}{x^6}+396 x^5+\frac{396}{x^5} -\frac{1521 x^4}{2}
 -\frac{1521}{2 x^4} +4082 x^3+\frac{4082}{x^3} -\frac{6541x^2}{2}-\frac{6541}{2 x^2} \\
 \hphantom{P_{2}(x) =}{}
 -8334 x-\frac{8334}{x}+16831.
\end{gather*}
Substituting them into \eqref{eq1.1.2} results in
\begin{gather*}	
2F\big(x,{\rm e}^{\hbar}\big) = \Big( x^{11/2} - \frac{1}{x^{11/2}}+ 2x^{15/2} - \frac{2}{x^{15/2}} + 5x^{19/2} - \frac{5}{x^{19/2}} + 13 x^{23/2} - \frac{13}{x^{23/2}}\\
\hphantom{2F\big(x,{\rm e}^{\hbar}\big) =}{}
+ 34 x^{27/2}- \frac{34}{x^{27/2}} - x^{29/2} + \frac{1}{x^{29/2}} + 89 x^{31/2} - \frac{89}{x^{31/2}} -2 x^{33/2} + \frac{2}{x^{33/2}}\\
\hphantom{2F\big(x,{\rm e}^{\hbar}\big) =}{}
+ 233 x^{35/2} - \frac{233}{x^{35/2}} -5x^{37/2} + \frac{5}{x^{37/2}} + \cdots \Big)\\
\hphantom{2F\big(x,{\rm e}^{\hbar}\big) =}{}
 + \hbar \Big( 5x^{11/2} - \frac{5}{x^{11/2}}+ 12x^{15/2} - \frac{12}{x^{15/2}} + 35x^{19/2} - \frac{35}{x^{19/2}} + 104 x^{23/2}\\
\hphantom{2F\big(x,{\rm e}^{\hbar}\big) =}{}
 - \frac{104}{x^{23/2}} + 306 x^{27/2} - \frac{306}{x^{27/2}} - 15x^{29/2} + \frac{15}{x^{29/2}} + 890 x^{31/2} - \frac{890}{x^{31/2}}\\
\hphantom{2F\big(x,{\rm e}^{\hbar}\big) =}{}
- 36 x^{33/2} + \frac{36}{x^{33/2}} + 2563 x^{35/2} - \frac{2563}{x^{35/2}} -105 x^{37/2} + \frac{105}{x^{37/2}} + \cdots \Big)\\
\hphantom{2F\big(x,{\rm e}^{\hbar}\big) =}{}
 + \hbar^2 \Big( \frac{25}{2} x^{11/2} - \frac{25}{2}\frac{1}{x^{11/2}}+ 36 x^{15/2} - \frac{36}{x^{15/2}} + \frac{247}{2} x^{19/2} - \frac{247}{2} \frac{1}{x^{19/2}}\\
\hphantom{2F\big(x,{\rm e}^{\hbar}\big) =}{}
+ 426 x^{23/2} - \frac{426}{x^{23/2}} + 1441 x^{27/2} - \frac{1441}{x^{27/2}} - \frac{225}{2} x^{29/2} + \frac{225}{2} \frac{1}{x^{29/2}} + 4781 x^{31/2}\\
\hphantom{2F\big(x,{\rm e}^{\hbar}\big) =}{}
- \frac{4781}{x^{31/2}}- 324 x^{33/2} + \frac{324}{x^{33/2}} + 2563 x^{35/2} - \frac{2563}{x^{35/2}} -\frac{2207}{2} x^{37/2} \\
\hphantom{2F\big(x,{\rm e}^{\hbar}\big) =}{}
+ \frac{2207}{2} \frac{1}{x^{37/2}} + \cdots \Big).
\end{gather*}	
Comparing to the series of the figure eight knot~\cite{GM}, we notice that every order of $\hbar$ appears in the above series whereas the series corresponding to the figure eight knot consists of only even powers of $\hbar$ (i.e., $P_{i}(x)=0$ for $i$ odd). This difference is an effect of the torus knot whose expansion involve all powers of $\hbar$~\cite{GM}. Furthermore, the $x$-terms begin from $m=11$ instead of $m=1$ and there are gaps in their powers. Specifically, $x^{\pm 13/2}$, $x^{\pm 17/2}$, $x^{\pm 21/2}$ and $x^{\pm 25/2}$ are absent. This is a consequence of the structure of \eqref{eq3.3.2}. A distinctive feature of the cable knot is that from~$x^{\pm 29/2}$ the coefficients are negative. Moreover, the positive and the negative coefficients alternate from that $x$-power for all $\hbar$ powers. These differences persist in the higher $\hbar$-orders. We will see these differences in a manifest way in the next section.

\section{Effects of the cabling}\label{sec6}

Since the initial data plays a core role in the recursion relation method, we discuss their features for the cable knot and then propose conjectures about it, which can be a useful guide for finding initial data for a family of the cable knots.

In the initial data (see \cite{C}) for the recursion relation \eqref{eq4.4.3}, we notice several differences from that of the figure eight knot~\cite{GM}. Before discussing them, let us begin with the properties of the~$F_K$ that are preserved by the cabling. The initial data consists of an odd number of terms and power of $q$ increases by one between every consecutive terms in a fixed $f_{m}$ for all $m$'s, which are also true for $f_{99}$ and $f_{101}$. Additionally, the reflection symmetry of coefficients is retained up to~$f_{43}$ for positive coefficients and up to~$f_{61}$ for the negative ones but of course, those $f_{m}$'s do not have the complete amphichiral structure. These invariant properties are a remnant feature of the amphichiral property of the figure eight knot.

A difference is that the nonzero initial data begins from $f_{11}$ and the gaps between the powers of~$x$ is four up to $x^{27/2}$, which is in the accordance with $f_{m}$'s. These features are direct consequences of the symmetric expansion of the Alexander polynomial of the cable knot \eqref{eq3.3.2}. In the case of the figure eight its coefficient functions start from $f_1$ and there are no such gaps. Another distinctive difference is that $f_{m}$'s containing negative coefficients appear from $m=29$. Moreover, the positive and negative coefficient $f_{m}$'s alternate from $f_{27}$ (i.e., positive coefficients for $f_{27}, f_{31}, \dots$ and negative coefficients for $f_{29}, f_{33},\dots$). Furthermore, from~$f_{47}$ the reflection symmetry of the positive coefficients in the appropriate $f_{m}$'s is broken. This phenomenon also occurs for the negative coefficient $f_{m}$'s from $m=65$. Breaking of the symmetry is expected since the cable knot of the figure eight is not amphichiral. The next difference is that the maximum power of $q$ in the positive coefficient~$f_{m}$'s for $m \geq 15$, the powers increase by $+2,+2,+3,+3,+4,+4,\dots, +11,+11$. For example, for $f_{15}$, $f_{19}$, $f_{23}$, $f_{27}$, and~$f_{31}$, their maximum powers are $q^6$, $q^{8}$, $q^{10}$, $q^{13}$ and $q^{16}$, respectively~\cite{C}. For the negative coefficient case, the changes in maximum powers are $4,4,5,5,6,6,\dots, 11,11$ from $m=33$. The minimum powers of~$f_m$'s having positive coefficients exhibit their changes as $0,0,-1,-1,-2,-2,-3,-3, \dots$ and for those with negative coefficients the pattern is $+2,+2,+1,+1,0,0,-1,-1,-2,-2,\dots$.

An universal feature of the negative coefficient $f_{m}$'s in the initial data is that their coefficient modulo sign is determined by the positive coefficient $f_{m}$. For example, the absolute value of the coefficients of $f_{29}$ is same as that of $f_{11}$; $f_{33}$'s coefficients come from that of $f_{15}$ up to sign and so forth. Hence coefficients of $f_{m}$ having negative coefficients are determined by $f_{m-18}$. In fact, this peculiar coefficient correlation also exists in the non-initial data $f_{101}$ whose coefficients are correlated with that of $f_{83}$.

\begin{Conjecture} For a class of a cable knot of the figure eight $K_{r}=C_{(r,2)}(4_1) \subset S^3$, ${r > 8}$ and odd, having monic Alexander polynomial, the coefficient functions $\big\{ f_{m}(q) \in \mathbb{Z} \big[q^{\pm 1}\big] \big\} $ of $F_{K_{r}}(x,q)$ can be classified into two (disjoint) subsets: one of them consists of elements having all positive coefficients $\big\{ f^{+}_{t}(q) \big\}_{t \in I^{+}}$ and the other subset contains elements whose coefficients are all negative $\big\{ f^{-}_{w}(q) \big\}_{w \in I^{-}}$. Furthermore, for every element in $\big\{ f^{-}_{w}(q) \big\}$, its coefficients coincide with that of an element in $\big\{ f^{+}_{t}(q) \big\}$ up to a sign.
\end{Conjecture}

\section{A relation to the figure eight knot}\label{sec7}

In this section, we observe an interesting relation between $F_{C_{(9,2)}(4_1)} (x,q)$ and $F_{4_1} (x,q)$. The latter was computed in~\cite{GM}. The relation enable us to circumvent the recursion method and hence provides an alternative and efficient method for computing $F_{C_{(9,2)}(4_1)} (x,q)$.

\begin{Proposition}
The positive and negative subsets of the initial data for \eqref{eq4.4.3} of $F_{C_{(9,2)}(4_1)} (x,q)$ are related to $F_{4_1} (x,q)$ as follows, respectively.
\begin{gather*}	
f_{11} (q) = h_1 (q) q^5,\\
f_{15} (q) = h_3 (q) q^6,\\
f_{19} (q) = h_5 (q) q^7,\\[-1mm]
\qquad\vdots\\
f_{43} (q) = h_{17} (q) q^{13},\\
f_{47} (q) = h_{19} (q) q^{14} + h_1 (q) q^{34}, \\
f_{51} (q) = h_{21} (q) q^{15} + h_3 (q) q^{39}, \\
f_{55} (q) = h_{23} (q) q^{16} + h_5 (q) q^{44}, \\[-1mm]
\qquad\vdots\\
f_{79} (q) = h_{35} (q) q^{22} + h_{17} (q) q^{74}, \\
f_{83} (q) = h_{37} (q) q^{23} + h_{19} (q) q^{79}+ h_1 (q) q^{99}, \\
f_{87} (q) = h_{39} (q) q^{24} + h_{21} (q) q^{84}+ h_3 (q) q^{108}, \\
f_{91} (q) = h_{41} (q) q^{25} + h_{23} (q) q^{89}+ h_5 (q) q^{117}, \\
f_{95} (q) = h_{43} (q) q^{26} + h_{25} (q) q^{94}+ h_7 (q) q^{126}, \\
f_{29} (q) = - h_1 (q) q^{15}, \\
f_{33} (q) = - h_3 (q) q^{18}, \\
f_{37} (q) = - h_5 (q) q^{21}, \\[-1mm]
\qquad\vdots\\
f_{61} (q) = - h_{17} (q) q^{39},\\
f_{65} (q) = - h_{19} (q) q^{42} - h_1 (q) q^{62}, \\
f_{69} (q) = - h_{21} (q) q^{45} - h_3 (q) q^{69}, \\
f_{73} (q) = - h_{23} (q) q^{48} - h_5 (q) q^{76}, \\[-1mm]
\qquad\vdots\\
f_{97} (q) = - h_{35} (q) q^{66} - h_{17} (q) q^{118},
\end{gather*}
where $h_{s}(q)$ are the coefficient functions of $F_{4_1}(x,q)$ $($see Appendix~{\rm \ref{appendixB})}.
\end{Proposition}

We note that $f_{1}=f_{3}=f_{5}=f_{7}=f_{9}=0$ as written in \cite{C}. We emphasize that the initial data for \eqref{eq4.4.3} in~\cite{C} were found from~\eqref{eq1.1.2}. The data turns out to be related to that of the figure eight knot. The above relation persists for $f_m (q)$ that are in the complement of the initial data set, namely, for $m > 97$. For example,
\begin{gather*}	
f_{99} (q) = h_{45} (q) q^{27} + h_{27} (q) q^{99}+ h_9 (q) q^{135}, \\
f_{101} (q) = - h_{37} (q) q^{69} - h_{19} (q) q^{125} - h_{1} (q) q^{145}.
\end{gather*}
They are in agreement with that obtained from \eqref{eq4.4.3}. We state the following conjectures.
\begin{Conjecture}
For a $(r,2)$-cabling of $4_1$ $( |r| \geq 2, \gcd (r,2)=1)$ in $\mathbb{Z} HS^3$, the coefficient functions $\{ f_m (q) \}$ of the series invariant $F_{C_{(r,2)}(4_1)}(x,q)$ are determined by the coefficient functions $\{ h_s (q)\} $ of $F_{4_1}(x,q)$ as
\begin{equation*}
f_m (q) = \pm \big( h_{s_1}(q) q^{w_1} + \cdots + h_{s_j}(q) q^{w_j} \big),\qquad w_i \in \mathbb{Z},
\end{equation*}
for some $j \in \mathbb{Z}_{+}$ $(j$ implicitly depends on $m)$.
\end{Conjecture}

More generally, we propose that
\begin{Conjecture}
For a $(r,2)$-cabling of any $($prime$)$ hyperbolic knot $K$ $(|r| \geq 2, \gcd(r,2)=1)$ in $\mathbb{Z} HS^3$, the coefficient functions $\{ f_m (q)\}$ of the series invariant $F_{C_{(r,2)}(K)}(x,q)$ are determined by the coefficient functions $\{ h_s (q)\}$ of $F_{K}(x,q)$ as{\samepage
\begin{equation*}
f_m (q) = \pm \big( h_{s_1}(q) q^{w_1} + \cdots + h_{s_j}(q) q^{w_j} \big),\qquad w_i \in \mathbb{Z},
\end{equation*}
for some $j \in \mathbb{Z}_{+}$.}
\end{Conjecture}

We finish by listing two applications of our result for future work. First, we can use $F_{C_{(p,2)}(4_1)}(x,q)$ to find $\hat{Z}$ associated with a closed oriented 3-manifold obtained by the Dehn surgery on the cable knot using the surgery formula in~\cite{GM, P}. This would in turn enable us to find the WRT invariant of the manifold using the result in~\cite{GPPV}. In both applications, it would extend $\hat{Z}$ and the WRT invariant to broader classes of 3-manifolds.

\appendix
\section{Appendix}

\subsection{The definitions of the operators}\label{appendixA}
We list the definitions of the operators in the $\hat{A}$-polynomial \eqref{eq4.4.1}:
\begin{gather*}
Q_2(t,M) = \hat{P}_2(t,t^4 M^2) \hat{P}_1\big(t,t^2 M^2\big) \hat{P}_0\big(t,t^6 M^2\big),\\
Q_1(t,M) = \hat{P}_0(t,t^4 M^2) \hat{P}_1\big(t,t^6 M^2\big) \hat{P}_2(t,t^2 M^2)- \hat{P}_1\big(t, t^6 M^2\big) \hat{P}_1\big(t,t^2 M^2\big) \hat{P}_1\big(t,t^4 M^2\big)\\
\hphantom{Q_1(t,M) = }{} + \hat{P}_2\big(t,t^4 M^2\big) \hat{P}_1\big(t,t^2 M^2\big) \hat{P}_0\big(t,t^6 M^2\big),\\
Q_0(t,M) = \hat{P}_0\big(t,t^4 M^2\big) \hat{P}_1\big(t, t^6 M^2\big) \hat{P}_0\big(t,t^2 M^2\big),
\\
\hat{P}_0(t,M):= t^6 M^4 \big({-}1 + t^{12} M^4\big),\\
\hat{P}_1(t,M):= -\big({-}1+ t^4 M^2\big)\big(1+ t^4 M^2\big) \big( 1-t^4 M^2 - t^4 M^4 - t^{12}M^4- t^{12}M^4 \\
\hphantom{\hat{P}_1(t,M):=}{}
- t^{12}M^6 + t^{16}M^8 \big),\\
\hat{P}_2(t,M):= t^{10}M^4 \big( {-}1 + t^4 M^4 \big).
\end{gather*}

\subsection{The data for the figure eight knot}\label{appendixB}

We record the initial data and the recursion relation for the figure eight knot from~\cite{GM}:
\begin{gather*}
h_1 (q) = 1,\\
h_3 (q) = 2,\\
h_5 (q) = \frac{1}{q}+3+q, \\
h_7 (q) = \frac{2}{q^2}+\frac{2}{q}+5+2 q+2 q^2, \\
h_9 (q) = \frac{1}{q^4}+\frac{3}{q^3}+\frac{4}{q^2}+\frac{5}{q}+8+5 q+4 q^2+3 q^3+q^4, \\
h_{11} (q) = \frac{2}{q^6}+\frac{2}{q^5}+\frac{6}{q^4}+\frac{7}{q^3}+\frac{10}{q^2}+\frac{10}{q}+15+10 q+10 q^2+7 q^3+6 q^4+2 q^5+2 q^6, \\
h_{13} (q) = \frac{1}{q^9}+\frac{3}{q^8}+\frac{4}{q^7}+\frac{7}{q^6}+\frac{11}{q^5}+\frac{15}{q^4}+\frac{18}{q^3}+\frac{21}{q^2}+\frac{23}{q}+27+23 q+21 q^2+18 q^3+15 q^4\\
\hphantom{h_{13} (q) =}{} +11 q^5+7 q^6+4 q^7+3 q^8+q^9,
\\
h_{m+14} (q) = -\frac{q^{-\frac{m}{2}-\frac{11}{2}}}{q^{\frac{m}{2}+\frac{13}{2}}-1} \big[ h_m \big(q^{\frac{m}{2}+\frac{17}{2}}-q^{m+9}\big)+ h_{m+2} \big(q^{\frac{m}{2}+\frac{15}{2}}-q^{\frac{m}{2}+\frac{17}{2}}+q^{m+9}-q^{m+10}\big) \\
\hphantom{h_{m+14} (q) =}{}
 + h_{m+4} \big({-}q^{\frac{m}{2}+\frac{11}{2}}-q^{\frac{m}{2}+\frac{17}{2}}-q^{\frac{m}{2}+\frac{19}{2}}+q^{\frac{3m}{2}+\frac{21}{2}}+q^{m+8}+q^{m+9}+q^{m+12}-q^7\big)\\
 \hphantom{h_{m+14} (q) =}{}
 + h_{m+6} \big({-}q^{\frac{m}{2}+\frac{9}{2}}+q^{\frac{m}{2}+\frac{11}{2}}-q^{\frac{m}{2}+\frac{15}{2}}-q^{\frac{m}{2}+\frac{17}{2}} +q^{\frac{3 m}{2}+\frac{25}{2}}+q^{m+9}+q^{m+10}\\
 \hphantom{h_{m+14} (q) ={}{}+ h_{m+6}}{}
 -q^{m+12}+q^{m+13}-q^5\big)\\
 \hphantom{h_{m+14} (q) =}{}
 + h_{m+8} \big(q^{\frac{m}{2}+\frac{11}{2}}+q^{\frac{m}{2}+\frac{13}{2}}-q^{\frac{m}{2}+\frac{17}{2}}+q^{\frac{m}{2}+\frac{19}{2}} -q^{\frac{3 m}{2}+\frac{31}{2}}-q^{m+8}+q^{m+9}\\
\hphantom{h_{m+14} (q) ={}{} + h_{m+8}}{}
 -q^{m+11}-q^{m+12}+q^2\big)\\
 \hphantom{h_{m+14} (q) =}{}
 + h_{m+10} \big(q^{\frac{m}{2}+\frac{9}{2}}+q^{\frac{m}{2}+\frac{11}{2}}+q^{\frac{m}{2}+\frac{17}{2}}-q^{\frac{3m}{2}+\frac{35}{2}}-q^{m+9}-q^{m+12}-q^{m+13}+1\big)\\
 \hphantom{h_{m+14} (q) =}{}
 + h_{m+12} \big(q^{\frac{m}{2}+\frac{11}{2}}-q^{\frac{m}{2}+\frac{13}{2}}+q^{m+11}-q^{m+12}\big)\big].
\end{gather*}

\subsection*{Acknowledgements} I would like to thank Sergei Gukov, Thang L\^e and Laura Starkston for helpful conversations. I~am grateful to Ciprian Manolescu for valuable suggestions on a draft of this paper. I~am also grateful to Colin Adams for valuable comments. I would like to thank to the referees for the suggestions that led to an improvement of my manuscript.


\pdfbookmark[1]{References}{ref}

\LastPageEnding

\end{document}